\input amstex
\documentstyle{amsppt}

\loadeufb
\loadeusb
\loadeufm
\loadeurb
\loadeusm

\magnification =\magstep 1
\refstyle{A}
\NoRunningHeads

\topmatter
\title On a problem of Piatetski-Shapiro and Shafarevich
\endtitle

\author Robert  Treger \endauthor
\address Princeton, New Jersey \endaddress
\email roberttreger117{\@}gmail.com \endemail

\endtopmatter

\document
\head  Introduction \endhead

The aim of this note is to revisit classical articles on the algebraic approach to uniformization and automorphic functions by Piatetski-Shapiro and Shafarevich. See \cite{S}\, (with comments and additional references), \cite{P}, and \cite{Pa-S, Sect.\, 4.1};  similar ideas were developed by Shimura independently of \cite{S}.

In classical geometric theory of automorphic forms, one
considers a bounded symmetric domain $V \subset {\bold C}^n$ and study automorphic
forms with respect to a discrete subgroup $\Gamma \subset Aut(V)$ in the group of all
complex analytic automorphisms of  $V$. Together with Siegel \cite{Si, vol.\, III, Sect.\, 6.1, 6.2}, one can ask about the scope of the classical theory. 

	From the geometric point of view, the short answer is a classification of the least complicated varieties of general type, namely, nonsingular projective algebraic varieties of general type with large and residually finite fundamental group as well as their noncompact counterparts.

Throughout the note, the ground field $k=\bold C$. After preliminaries, we consider several important examples. 

In the last section, we prove our main theorem that provides a partial solution of a problem proposed by  by Piatetski-Shapiro and Shafarevich in \cite{S, Introduction}, \cite{P}, and \cite{Pa-S, p.\, 83}, namely, whether the existence of proalgebraic quasi-homogeneous coverings of general type is the characteristic property of algebraic varieties whose universal coverings are bounded symmetric domains.

\head 1. Preliminaries
\endhead

\subhead 1.1  \endsubhead
All proalgebraic varieties will be assumed {\it finite-dimensional}, {\it irreducible}, and {\it normal}. Recall that a
proalgebraic variety can be represented as a projective limit
$$
Y=\varprojlim X_\alpha  
$$
of a countable set of algebraic varieties $X_\alpha$, where $f_{\alpha,\beta}:
X_\alpha \to X_\beta$ are integral morphisms \cite{S, Sect.\, 4, Proposition 2}. A point
$x\in Y$ is called {\it interior}\/ if the corresponding maximal ideal $m_x$ has a
finite number of generators or, equivalently, there exists an index $\gamma$ such
that all morphisms $X_\alpha \to X_\gamma$ are unramified at $x_\gamma$, where
$x_\gamma$ is the projection of $x$ to $X_\gamma$ \cite{S, Sect. 4, Definition 4, Proposition 3}. We also call a point $x_\alpha\in X_\alpha$ interior if
the point $x\in Y$ lying over it is interior. 

	Unless stated otherwise, throughout the note  $Y$ is {\it not}\/  an {\it algebraic}\/ variety. 

\subhead 1.2  \endsubhead
We assume the set of interior points of $Y$, denoted by $U_Y$, is not empty, nonsingular, and $\pi^{alg}_1(U_Y)=\{1\}$. Then there exists an index $\gamma$
such that for $\alpha \geq \gamma$ all interior points of $X_\alpha$, denoted by  $U_\alpha$,
are nonsingular. In particular, $U_\alpha$ is a complex manifold. 

In the sequel, we also assume $\pi_1(U_\alpha)$ is {\it residually finite} for every index $\alpha$, and {\it large} \cite{Kol}. 
The scheme $U_Y$ has another topology, so-called pro-etale topology, so that the projections $U_Y \to U_\alpha$ are continuous in the etale topology on $U_\alpha$. Let $\tilde U$ denote the universal topological covering of $U_\alpha$. Then $\tilde U$ is also a complex manifold.

 We get a natural topological embedding
$\tilde U \hookrightarrow U_Y$, where $U_Y$ is equipped with its pro-etale topology. Our $\tilde U$ is also equipped with a Zariski-type topology induced by the Zariski topology on $U_Y$.

\subhead 1.3 \endsubhead Let $K$ be a field of finite transcendence degree
$n$ over
$k$.  We denote by $\bold D ^\ell(K)$ the set of all regular
$n$-differentials of $K$ of weight $\ell$ \cite{S, Sect.\, 5}. Let $\bold D= \sum \bold
D^\ell(K)$ be the graded algebra of differentials. If $\bold D$ is irreducible and
integrally closed, and all its elements are integral over a certain homogeneous
subalgebra of finite type then $Y=\text{Proj}\, \bold D$ is a {\it projective proalgebraic} variety \cite{S, Sect.\, 4, Definition 3}. 

Let $k(Y)$ denote the field of rational function on $Y$. Clearly always $k(Y) \subseteq K$. Let $\eusm K_\alpha$ denote the canonical bundle
on $ U_\alpha \subset X_\alpha$, for a cofinal system of indexes $\alpha$. 

\definition {Definition 1.3.1} The field  $K$ is called a field of general type (in the sense of Piatetski-Shapiro and Shafarevich) if $k(Y)=K$ and there is an integer $m>0$ such that each $\eusm K^m_\alpha$ and its global sections define an embedding in the corresponding 
finite-dimensional projective space. 
\enddefinition

\head 2. Examples 
\endhead
 
\example {Example 1}   Let $V \subset{\bold C}^n$  be a bounded domain, and $\Gamma$ be a fixed point free discrete subgroup of 
$ Aut(V)$ such that $V/\Gamma$ is a compact complex manifold.
Let $k(\Gamma) = \cup K_{\Delta_i}$, where $\Delta_i$ runs over 
subgroups of finite index in  $\Gamma$ ($\cap\Delta_i = 1$).  Employing Poincar\'e series, it was shown that each $V/\Delta_i $ is a nonsingular projective variety with ample canonical bundle, and it is the {\it absolute minimal}\/ model (Shioda) of its field of rational functions. Therefore
$$
\varprojlim   V_\alpha/\Delta_i \simeq \text{Proj}\, \bold D, \qquad \bold D=\sum \bold
D^\ell (k(\Gamma)),
$$
and all points of $\text{Proj}\, \bold D$ are interior.  
\endexample

\example {Example 2(i)} Let $V \subset{\bold C}^n$  be a
bounded symmetric domain and $\Gamma \subset Aut(V)$ be a fixed point free arithmetic subgroup
in the group of all complex analytic automorphisms of $V$. Let $k(\Gamma) = \cup K_{\Delta_\alpha}$, where $\Delta_\alpha$ runs over neat
subgroups of finite index in  $\Gamma$ ($ \cap\Delta_\alpha = 1$). Then each 
$V/\Delta_\alpha$ is a nonsingular subvariety of its Baily-Borel compactification (see, e.g., [M, Prop 3.4]).

According to Tai, each $V/\Delta_\alpha$ is a variety of general type provided
$\Delta_\alpha$ is sufficiently small. The field $k(\Gamma)$ has abundance of regular $n$-differentials. Namely, the regular 
$n$-differentials on
$V/\Delta_\alpha$ that have at worst logarithmic poles along the boundary of $V/\Delta_\alpha$
generate the homogeneous coordinate ring of the Baily
-Borel compactification $\overline{V/\Delta_\alpha}$ \cite{M, Prop 4.2}. Further, it follows  from Mumford's proof of Tai's theorem that these differentials are regular differentials of $k(\Gamma)$ because each $\overline{V/\Delta_\alpha}$ has a tower of
coverings {\it universally ramified} over its boundary (for details, see \cite{M, Sect.\, 4}).
\endexample

\example{Example 2(ii)}  Let $V \subset{\bold C}^{3g-3}$, $g\geq2$,
be the Teichmuller moduli space of curves of genus $g$. Let $\Gamma \subset
\Gamma_g$ be a sufficiently small subgroup of the mapping class group
$\Gamma_g$ that acts freely on $U$. This example is similar to Example 2.2(i) as was suggested in \cite{M, Sect.\, 4}.  One can apply a more recent
result of Looijenga \cite{L} who finally established that the mapping class group has the key property, namely, {\it local universal ramification over the boundary},\/  needed to show that $k(\Gamma)$ has abundance of regular $n$-differentials as in Example 2(i). 
\endexample

\example {Example 3} Let $X$ be a nonsingular $n$-dimensional projective variety with {\it ample} canonical bundle $\Cal K_X$, and {\it large} and {\it residually finite} fundamental group, and without general elliptic curvilinear sections. 
Assume the universal covering $V$ of $X$ has a $q$-{\it Bergman}\/ metric for an integer $q$. Hence $V$ is a bounded domain \cite{T3, Corollary}. Thus we are in the situation of Example 1.
\endexample

\example {Example 4} (See Remark 4.3 by Campana in \cite{T1}.) Let $X$ be a sufficiently ample divisor in an Abelian variety of dimension at least 3. The universal covering $V$ of $X$ is not a bounded domain, and the fundamental group $\pi_1(X)$ is an  Abelian group of rank at least 3, in particular, amenable. The canonical bundle on $X$ is ample.
If dim$A=3$ then $Aut(V)^\circ=\{1\}$ according to Nadel \cite{N}. 

Since $V$ is Stein, $X$ is an absolute minimal model of its field of rational functions according to [Kob, (6.3.21)] (generalizing earlier partial results by Igusa and Shioda).
\endexample

\head
3. A characterization of algebraic varieties whose universal coverings are bounded symmetric domains
\endhead

\subhead 3.1\endsubhead The problem \lq\lq whether the existence of quasi-homogeneous proalgebraic coverings of general type (in the sense of Piatetski-Shapiro and Shafarevich) is the characteristic property of algebraic varieties whose universal coverings are  bounded symmetric domains\rq\rq\/ is stated in \cite{S, Introduction}, \cite{P},  \cite{Pa-S, pp.\, 82-83}. Recall 

\definition{Definition 3.1.1  \cite{S, Sect.\, 6}} A proalgebraic variety $Y$ is called
{\it quasi-homogeneous}\/  if the set of its boundary points is different from $Y$ and is closed, and the orbit of every interior point relative to the group $Aut(Y)$ of all automorphisms of $Y$ is everywhere dense (in Zariski topology) in $Y$.
\enddefinition
If $Y$ is quasi-homogeneous then the set of its nonsingular
points coincides with the set of interior points \cite{S, Sect.\, 6, Proposition 2}.

\subhead 3.2 
\endsubhead
 Let $Y=\varprojlim X_\alpha$ be a projective proalgebraic variety. We keep the notation and assumptions of the preliminaries.

 Let
$k(Y)=\cup k(X_\alpha)$ where $k(Y)$ and $k(X_\alpha)$ are the fields of rational functions on $Y$ and $X_\alpha$, respectively. We assume $Y=\text{Proj}\,\sum \bold D^q(k(Y))$. Recall that $k(X_\alpha)$ is the field of
rational functions on $U_\alpha =\tilde U/\Delta_\alpha$ for all $\alpha$ greater than
a fixed sufficiently large index $\gamma$, where $\Delta_\alpha$ is a subgroup in the group $Aut(\tilde U)$ of complex analytic automorphisms of $\tilde U$.  Let 
$$
Comm(Y) : = Comm(\Delta_\gamma) \subset Aut(\tilde U)
$$ 
denote the \it commensurability\/ \rm subgroup.

 Let $G$ denote  the identity component of the closure of  $ Comm(Y)$ in $Aut(\tilde U)$. Given a point $x\in \tilde U$, let $G/B_x$ denote the orbit of $G$ through $x$ where $B_x \subset G$ is the stabilizer of  $x$. 
The definition as well as non-vanishing of the cohomology group $H^{top}(G/B_x)$ is discussed in \cite{H, p.\, 890} and \cite{Kos, Theorem 13.1}. We observe that $G$ acts effectively on $G/B_x$, and $B_x$ is reductive in $G$ (see \cite{H, (1.1)}).

\subhead 3.3
\endsubhead 
What follows is a global version of \cite{FK, pp.\, 5-15, Proposition I.1.6, Proposition I.2.1}. We  replaced functions by sheaves on manifolds. Classically, embeddings in infinite-dimensional projective spaces  were considered by Bochner, Calabi \cite{C}, and in articles by Kobayashi (see, e.g., \cite{Kob, Chap.\, 4.10}).

We keep the notation and assumptions of (1.1) - (1.3). Assume $\eusm K^m_\alpha$ and its global sections define the embedding of $U_\alpha$ in the corresponding finite-dimensional projective space, and $k(Y)$ is of general type (as in Definition 1.3.1 with the same $m$). As in \cite{T1, T2}, we get the corresponding  {\it section} 
$$
\bold B:=B_{\tilde U, \eusm K^m}(z,w)
$$
 holomorphic in $z$ and antiholomorphic in $w$, where $\eusm K$ is the canonical bundle on $\tilde U$. Further, $\bold B(z,z)>0$ and $\log \bold B(z, z)$ is strictly plurisubharmonic. A priory, $\bold B$ is {\it not}\/ a Hermitian kernel  of {\it positive type} \cite{FK, p.\,8, p.\, 12}.

We denote by $\bold B^o$ the restriction of $\bold B$ to the orbit $G/B_x$. Clearly, a Hilbert subspace of a Hilbert space with a reproducing kernel also has a reproducing kernel.  So $\bold B^o$ is a {\it reproducing kernel}\/ of a unique Hilbert space $H$ of holomorphic sections. Indeed, that is true locally in a Euclidean neighborhood of any point of $G/B_x$ hence globally because $G/B_x$ is homogeneous.

We consider the projective space $\bold P(H^*)$ with its Fubini-Study metric. In local coordinates, the evaluation at a point $Q\in G/B_x$,
$
e_Q:  f\mapsto f(Q),
$
is a continuous linear functional on $H$.
 We get a natural map 
$$
\Upsilon : G/B_x \longrightarrow \bold {P}(H^*).
$$
By the assumption, $\Upsilon$ is an embedding. Let $Aut(G/B_x)$ denote the group of complex analytic automorphisms of $G/B_x$. All elements of  $Aut(G/B_x)$ act as collineations of $\bold P(H^*)$ (compare \cite {FK, p.\, 15}).

\remark {Remark 3.4} We are in the situation of the standard Bergman metric. Hence $G/B_x$ \lq\lq appears\rq\rq to be a bounded domain. On the other hand, sections like $\bold B$, above, exist in a more general situation because locally in Eucledian neighborhoods we have convergence like in \cite{T2}. However, we get a metric that may not be the Bergman metric.
\endremark
\proclaim{Theorem 3.5} Let $Y$ be a projective proalgebraic  variety of general type. We keep the notation and assumptions of $(1.1)-(1.3)$ and $(3.2)-(3.3)$.  Let $x\in \tilde U$ be a point such that the orbit of $G$ through $x$ is not zero-dimensional and is Zariski dense in $\tilde U$.  Then $\tilde U$ is a bounded symmetric domain.
\endproclaim

\demo{Proof} Let $G/B_x$ denote the orbit through $x\in\tilde U$. Then it has a complex homogeneous structure induced by its embedding in $\bold P(H^*)$. We will apply classical results of Borel, Koszul, and Hano.

By assumptions, this complex structure is Kahler homogeneous. Further, the Ricci form as well as Koszul's canonical form are non-degenerate (see \cite{H, Sect.\, 6, p.\, 893}), and $H^{top}(G/B_x)\not=0$. The latter follows from \cite{H, p.\, 890} and \cite{Kos, Theorem 13.1}  because of our embedding  $G/B_x \subset \bold P(H^*)$. Clearly, $G$ acts effectively on $G/B_x$.
 Hence we may apply a theorem of Hano \cite{H, Theorem 3 and p.\, 893}.  It follows $G$ is a semisimple group. 

 According to Borel \cite{B, Theorem 4}, $G/B_x$ is simply connected.  Furthermore, $G/B_x$ is a bounded symmetric domain since the corresponding fundamental group is large, i.e., $G/B_x$ contains no compact analytic subsets.

For an appropriate discrete subgroup $\Delta_\alpha$,  $\Delta_\alpha   \backslash G/B_x$ is a locally symmetric space and we may consider its Satake compactification (Satake, Baily, Borel, Piatetski-Shapiro).  It follows that its closure in $X_\alpha$ is an {\it algebraic} subvariety by a theorem of Chow.

Since the orbit is Zariski dense in $\tilde U$, we get $\tilde U$ is a bounded symmetric domain.
\enddemo 

%\remark {Remark 3.6} We conjecture that with the {\it basic} assumptions of the %theorem, the orbits can not be zero-dimensional and Zariski dense.
%\endremark

\enddocument